\documentclass{amsart}
\usepackage{graphicx}
\usepackage{amscd}
\usepackage{amsmath}
\usepackage{amsfonts}
\usepackage{amssymb}

\theoremstyle{plain}

\newtheorem{definition}{Definition}

\newtheorem{lemma}{Lemma}

\newtheorem{proposition}{Proposition}

\numberwithin{equation}{section}

\begin{document}
\title{A duality theory for some non-convex functions of matrices.}
\author{Ivar Ekeland, Canada Research Chair in Mathematical Economics}
\address{Department of Mathematics, University of British Columbia}
\email{ekeland@math.ubc.ca}
\thanks{The author thanks Jey Sivaloganathan for discussing an early version of the
paper, Paolo Marcellini for suggesting Example 3, and the participants of the
Pacific Northwest Optimization Workshop, May 2003}
\date{May 15, 2003}

\begin{abstract}
We study a special class of non-convex functions which appear in nonlinear
elasticity, and we prove that they have a well-defined Legendre transform.
Several examples are given, and an application to a nonlinear eigenvalue problem
\end{abstract}\maketitle

\section{Introduction}

We want to define a Legendre transform $F^{L}\left(  y\right)  $ for functions
$F\left(  x\right)  $, where $x$ is an $N\times K$ matrix, and $F$ involves
the various cofactors of $x$. Note that $F$ then has to be strongly nonlinear,
and nonconvex. The simplest case is when $N=K$ and $F\;$is a function of the
determinant only:\ $F\left(  x\right)  =\Phi\left(  \det x\right)  $. We show
that the Legendre transform of $\ $%
\[
F\left(  x\right)  =\frac{N}{p}\left|  \det x\right|  ^{p/N}%
\]
is
\[
F^{L}\left(  y\right)  =\frac{N}{q}\left|  \det y\right|  ^{q/N}%
\]
with $1/p+1/q=1$, thereby generalizing the classical duality between $L^{p}$
spaces. The next simplest is when $F$ is a function of the $\left(
N-1\right)  $-cofactors of $x$: in the case when $N=K=3$, we give conditions
under which $F$ has a well-defined Legendre transform. This covers for
instance the area functional, defined over $2\times3$ matrices $x=\left(
x_{j}^{i}\right)  $ by:
\[
F\left(  x\right)  =\left[  \left(  x_{1}^{1}x_{2}^{2}-x_{1}^{2}x_{2}%
^{1}\right)  ^{2}+\left(  x_{1}^{1}x_{3}^{2}-x_{1}^{2}x_{3}^{1}\right)
^{2}+\left(  x_{2}^{1}x_{3}^{2}-x_{2}^{2}x_{3}^{1}\right)  ^{2}\right]  ^{1/2}%
\]
which turns out to be its own Legendre transform (in other words, it is self-dual).

Note that with our definition, the Legendre transform $F^{L}$ of $F$ satisfies
the usual duality relations:
\begin{align*}
\left(  F^{L}\right)  ^{L}  & =F\\
\left(  F^{\prime}\right)  ^{-1}  & =\left(  F^{L}\right)  ^{\prime}%
\end{align*}

The paper is organized as follows. First we define precisely what we mean by a
Legendre transform. Then we sudy two polar cases. In the first one $F\;$is a
function of the determinant only, and in the second $F$ depends only on the
$2\times2$ cofactors of a $3\times3$ matrix. We give seeveral examples, and we
conclude by giving an application to a nonlinear eigenvalue problem.

\section{The Legendre transform.}

Let $X$ be a finite-dimensional vector space, and $Y$ its dual, the duality
pairing being denoted by $<x,y>$. Let $F:X\rightarrow R$ be a $C^{1}$ function
and $F^{\prime}\left(  x\right)  \in Y$ its derivative at $x$. The classical
formula of Legendre associates with every $y\in Y$ a set $\Gamma_{L}\left(
y\right)  \subset R$ defined as follows:
\begin{equation}
\Gamma_{L}\left(  y\right)  =\left\{  <x,y>-F\left(  x\right)  \mid
y=F^{\prime}\left(  x\right)  \right\} \label{Legendre}%
\end{equation}

Usually, the right-hand contains several points, so that $\Gamma_{L}$ is a
multi-valued map from $Y$ to $R$. We refer to \cite{Ivar1} for a study of this
map. For certains classes of functions $F$, however, the right-hand side is a
singleton, so that formula (\ref{Legendre}) defines a function on $Y$, which
is then called the Legendre transform of $F$.

\begin{definition}
Consider a map $F:$ $\Omega\rightarrow R$, where $\Omega$ is a submanifold of
$X$, and set $\Sigma=F^{\prime}\left(  \Omega\right)  $. We shall say that a
$C^{1}$ function $G:\Sigma\rightarrow R$ \ is the Legendre transform of $F$ if
$\Sigma$ is a submanifold of $Y$ and:
\[
\left[  x\in\Omega,\;y=F^{\prime}\left(  x\right)  \right]  \Longrightarrow
\;<x,y>-F\left(  x\right)  =G\left(  y\right)  \text{\ }%
\]
\end{definition}

Of course, if $F^{\prime}$ is one-to-one, this formula becomes:
\[
G\left(  y\right)  =<\left(  F^{\prime}\right)  ^{-1}\left(  y\right)
,y>-F\left(  \left(  F^{\prime}\right)  ^{-1}\left(  y\right)  \right)  \text{
}\forall y\in\left(  F^{\prime}\right)  ^{-1}\left(  \Omega\right)
\]
We shall denote the Legendre transform of $F$ by $F^{L},$ so that $F^{L}=G$ in
the above. It follows from the general theory of the Legendre transform (see
\cite{Ivar1}) that $F^{\prime}$and $G^{\prime}$ are inverse of each other and
that $F\left(  x\right)  $ is the Legendre transform of $F^{L}\left(
y\right)  $. In other words, we have the classic formulas:
\begin{align*}
\left(  F^{L}\right)  ^{L}  & =F\\
\left(  F^{\prime}\right)  ^{-1}  & =\left(  F^{L}\right)  ^{\prime}%
\end{align*}

There is a well-defined theory of Legendre transform for convex functions. In
that case, $\Omega=X$ and the Legendre formula is replaced by the Fenchel
formula:
\[
F^{\ast}\left(  y\right)  =\sup_{x}\left\{  yx-F\left(  x\right)  \right\}
\]
so that differentiability is no longer required. The Fenchel transform
$F^{\ast}$ will coincide with the Legendre transform $F^{L}$ provided $F$ is
$C^{1}$ and strictly convex. We now proceed to give other classes of functions
which have a well-defined Legendre transform. Roughly speaking, these will be
functions $F\left(  x\right)  $, where $x$ is a matrix and $F$ depends only on
the cofactors of $x$.

Given a number $K$ and some $k\leq K$, we shall denote by $\mathcal{P}_{K}$
the set of strictly increasing maps of $\left\{  1,...,k\right\}  $ into
$\left\{  1,...,K\right\}  $:
\[
\mathcal{P}_{K}=\left\{  \pi:\left\{  1,...,k\right\}  \rightarrow\left\{
1,...,K\right\}  \mid\pi\left(  1\right)  <...<\pi\left(  k\right)  \right\}
\]
and by $c\left(  K,k\right)  =C_{K}^{k}$ its cardinal. We can also think of
$\mathcal{P}_{K}$ as being the set of ordered subsets of $\left\{
1,...,K\right\}  $ with $k$ elements. Similar notations will hold for
$\mathcal{P}_{N}$ and $c\left(  k,N\right)  $, provided $k\leq N$.

Consider a fixed $N$-dimensional space $E.$ An element $x\in E^{\ast}$ will
have coordinates $x_{n}$, $1\leq n\leq N$. Given a family of $K$ linear forms,
$\left(  x^{1},...,x^{K}\right)  $, and a number $k\leq\min\left\{
K,N\right\}  $, there are $c\left(  K,k\right)  \times c\left(  N,k\right)  $
square $k\times k$ matrices which can be extracted from the matrix $x_{n}%
^{k}.$ Each of them is specified by a certain choice of $k$ lines and $k $
columns, that is by some $\pi\in\mathcal{P}_{K}$ and some $\sigma
\in\mathcal{P}_{K}$. We shall denote it by:
\[
x_{\sigma}^{\pi}=\left[  x_{\sigma\left(  j\right)  }^{\pi\left(  i\right)
}\right]  _{1\leq j\leq k}^{1\leq i\leq k}%
\]
and we shall denote by $\Delta_{k}^{\pi,\sigma}(x^{1},...,x^{K})$ its
determinant:
\[
\Delta_{k}^{\pi,\sigma}(x^{1},...,x^{K})=\det\left[  x_{\sigma}^{\pi}\right]
\]

We then define a map $\Delta_{k}:\left(  E^{\ast}\right)  ^{K}\rightarrow
R^{c\left(  K,k\right)  }\times$ $R^{c\left(  N,k\right)  }$ by:
\begin{equation}
\Delta_{k}=\left(  \Delta_{k}^{\pi,\sigma}\right)  _{\sigma\in\mathcal{P}_{N}%
}^{\pi\in\mathcal{P}_{K}}\label{e1}%
\end{equation}

So the map $\Delta_{k}$ just associates with a $N\times K$ matrix $\left(
x_{n}^{k}\right)  $ the determinants of all the $k\times k\,$\ matrices which
can be extracted from it, that is, its $k$-cofactors

\begin{lemma}
We have:
\[
\Delta_{k}^{\pi,\sigma}=k\sum_{n,j}x_{n}^{j}\frac{\partial\Delta_{k}%
^{\pi,\sigma}}{\partial x_{n}^{j}}%
\]
\end{lemma}

\begin{proof}
This is just the Euler identity for $k$-homogeneous functions
\end{proof}

\begin{lemma}
\label{lem1}Set:
\[
z\left(  \pi,\sigma\right)  _{j}^{n}=\frac{\partial}{\partial x_{j}^{n}}%
\det\left[  x_{\sigma}^{\pi}\right]
\]
We then have:
\[
\det\left[  z\left(  \pi,\sigma\right)  _{\sigma}^{\pi}\right]  =\left(
\det\left[  x_{\sigma}^{\pi}\right]  \right)  ^{k-1}%
\]
\end{lemma}

\begin{proof}
We know that $\left(  z_{\sigma}^{\pi}\right)  _{j}^{n}$ is $0$ if $n$ does
not belong to the image of $\pi$, or if $j$ does not belong to the image of
$\sigma$, and that otherwise it is just the cofactor of $x_{j}^{n}$ in the
matrix $\left(  x_{\sigma}^{\pi}\right)  $. The last identity then follows
from the well-known fact that the determinant of a $k\times k$ matrix raised
to the $\left(  k-1\right)  $-th power is the determinant of the\ cofactor matrix.
\end{proof}

As we stated in the beginning, we are interested in functions of the $\left(
x_{n}^{k}\right)  $ which involve the $\Delta_{k}$. We now make this idea
precise. Consider the map:
\begin{align*}
\Delta & =\left(  \Delta_{1},...,\Delta_{K}\right)  :R^{NK}\rightarrow H\\
H  & =R^{NK}\times...\times R^{c\left(  N,k\right)  c(K,k)}\times...\times
R^{c\left(  N,K\right)  }%
\end{align*}

A function $\Phi:H\rightarrow R$ will be called \emph{trivial} if it depends
on the $NK$ first coordinates only that is, if it factors through $\left(
E^{\ast}\right)  ^{K}$. Note that every function $F:R^{NK}\rightarrow R$ can
be written $F=\Phi\circ\Delta$, where $\Phi:H\rightarrow R\ $is the identity
on $\left(  E^{\ast}\right)  ^{K}$ and sends all the other coordinates to $0$.
This is called the trivial factorisation.

\begin{definition}
A function $F:R^{NK}\rightarrow R$ will be called $k$-\emph{adapted} if it
factors through $\Delta_{k}$, that is, if we have $F=\Phi\circ\Delta_{k}$ for
some function $\Phi:R^{c\left(  N,k\right)  c(K,k)}\rightarrow R$. It is
\emph{adapted} if it factors non-trivially through $\Delta$, that is, if we
have $F=\Phi\circ\Delta$ for some non-trivial function $\Phi:H\rightarrow R$
such that $\Phi\circ\Delta_{k}$
\end{definition}

We now write the formula for the Legendre transform. Take an adapted function
$F:R^{NK}\rightarrow R$:
\[
F\left(  x\right)  =F\left(  x^{1},...,x^{K}\right)  =F\left(  x_{n}%
^{j}\right)  =\Phi\left(  \Delta_{k}^{\pi,\sigma}\right)  =\Phi\left(
\Delta\right)
\]
and pair $K\times N$ matrices with $N\times K$ matrices by:
\[
<x,y>=\sum_{n,j}x_{n}^{j}y_{j}^{n}%
\]

Substitute in the definition (\ref{Legendre}):%

\begin{align}
\Gamma_{L}\left(  y_{1},...,y_{K}\right)   &  =\left\{  \sum_{n,j}x_{n}%
^{j}y_{j}^{n}-F\left(  x^{1},...,x^{K}\right)  \mid y_{j}^{n}=\frac{\partial
F}{\partial x_{n}^{j}}\right\} \nonumber\\
&  =\left\{  \sum_{n,j}x_{n}^{j}y_{j}^{n}-\Phi\left(  \Delta\right)  \mid
y_{j}^{n}=\sum_{k,\pi,\sigma}\frac{\partial\Phi}{\partial\Delta_{k}%
^{\pi,\sigma}}\frac{\partial\Delta_{k}^{\pi,\sigma}}{\partial x_{n}^{j}%
}\right\} \nonumber\\
&  =\left\{  \sum_{n,j,k,\pi,\sigma}x_{n}^{j}\frac{\partial\Phi}%
{\partial\Delta_{k}^{\pi,\sigma}}\frac{\partial\Delta_{k}^{\pi,\sigma}%
}{\partial x_{n}^{j}}-\;\Phi\left(  \Delta\right)  \mid y_{j}^{n}=\sum
_{k,\pi,\sigma}\frac{\partial\Phi}{\partial\Delta_{k}^{\pi,\sigma}}%
\frac{\partial\Delta_{k}^{\pi,\sigma}}{\partial x_{n}^{j}}\right\} \nonumber\\
&  =\left\{  \sum_{k,\pi,\sigma}\frac{\partial\Phi}{\partial\Delta_{k}%
^{\pi,\sigma}}\sum_{n,j}x_{n}^{j}\frac{\partial\Delta_{k}^{\pi,\sigma}%
}{\partial x_{n}^{j}}-\;\Phi\left(  \Delta\right)  \mid y_{j}^{n}=\sum
_{k,\pi,\sigma}\frac{\partial\Phi}{\partial\Delta_{k}^{\pi,\sigma}}%
\frac{\partial\Delta_{k}^{\pi,\sigma}}{\partial x_{n}^{j}}\right\} \nonumber\\
&  =\left\{  \sum_{k,\pi,\sigma}k\frac{\partial\Phi}{\partial\Delta_{k}%
^{\pi,\sigma}}\Delta_{k}^{\pi,\sigma}-\Phi\left(  \Delta\right)  \mid
y_{j}^{n}=\sum_{k,\pi,\sigma}\frac{\partial\Phi}{\partial\Delta_{k}%
^{\pi,\sigma}}\frac{\partial\Delta_{k}^{\pi,\sigma}}{\partial x_{n}^{j}%
}\right\} \label{m5}%
\end{align}
where $\Delta$ stands for $\Delta\left(  x^{1},...,x^{K}\right)  $. We rewrite
the result in more compact notation:
\begin{equation}
\Gamma_{L}\left(  y_{1},...,y_{K}\right)  =\left\{  \sum_{k}k\frac
{\partial\Phi}{\partial\Delta_{k}}\Delta_{k}-\Phi\left(  \Delta\right)  \mid
y_{j}^{n}=\sum_{k}\frac{\partial\Phi}{\partial\Delta_{k}}\frac{\partial
\Delta_{k}}{\partial x_{n}^{j}}\right\} \label{m1}%
\end{equation}

We would like to give general conditions on $\Phi$ which would ensure that the
right-hand side is a singleton, so that $F$ has a well-defined Legendre
transform. In addition, we would like to show that if $F$ is $k$-adapted, then
$F^{L}$ is $k$-adapted as well. Unfortunately, we have not been able to fulfil
this program (the calculations very quickly become horrendous) so we will be
content with two examples$.$

\section{Functions of the determinant.}

We take $N=K$. So let $x$ be the square matrix with coefficients $x_{n}^{k}$.
Denote by $X_{n}^{k}$ the cofactor of $x_{n}^{k}$ in $X$. We consider
functions $F:R^{N}\rightarrow R$ of the following type:%

\[
F\left(  x\right)  =\Phi\left(  \det x\right)
\]
where $\Phi:R\rightarrow R$ is a $C^{1}$ function.

Let us apply the preceding theory. We have:
\begin{equation}
y_{k}^{n}=\frac{\partial F}{\partial x_{n}^{k}}=\Phi^{\prime}\left(  \det
x\right)  \frac{\partial\det x}{\partial x_{n}^{k}}=\Phi^{\prime}\left(  \det
x\right)  \det X_{n}^{k}\label{li7}%
\end{equation}

Hence:
\begin{align}
\sum_{n,k}y_{k}^{n}x_{n}^{k}  & =\Phi^{\prime}\left(  \det x\right)
\sum_{n,k}x_{n}^{k}\det X_{n}^{k}\label{li1}\\
& =N\;\Phi^{\prime}\left(  \det X\right)  \det X\;\label{li2}%
\end{align}

On the other hand, denoting by $Y$ the matrix with coefficients $y_{k}^{n}$,
and by $z$ the matrix with coefficients $z_{k}^{n}=\det X_{n}^{k}$, we have:%

\begin{equation}
\det Y=\left(  \Phi^{\prime}\left(  \det x\right)  \right)  ^{N}\det
z=\Phi^{\prime}\left(  \det x\right)  ^{N}\left(  \det x\right)
^{N-1}\label{li4}%
\end{equation}

\begin{proposition}
Assume $\Phi$ is such that the function $t\rightarrow t^{N-1}\Phi^{\prime
}\left(  t\right)  ^{N}$ is invertible on $\left(  a,b\right)  $, and let
$\psi:\left(  a^{\prime},b^{\prime}\right)  \rightarrow R$ be its inverse. Set
$\Omega=\left\{  x\;|\;a<\det x<b\right\}  $ and $\Sigma=\left\{
y\;|\;a^{\prime}<\det y<b^{\prime}\right\}  $. Then the function
$F:\Omega\rightarrow R$ given by:
\[
F\left(  x\right)  =\Phi\left(  \det x\right)
\]
has a Legendre transform $F^{L}:\Sigma\rightarrow R$ given by:
\begin{equation}
F^{L}(y)=N\;\Phi^{\prime}\left(  \psi\left(  \det y\right)  \right)
\psi\left(  \det y\right)  \;-\Phi\left(  \psi\left(  \det y\right)  \right)
\label{li11}%
\end{equation}
\end{proposition}

\begin{proof}
Equation (\ref{li4}) gives $\det x=\psi\left(  \det y\right)  $. Writing
(\ref{li2}) and (\ref{li4}) back into formula (\ref{m5}), we get:
\begin{align*}
\Gamma_{L}\left(  y\right)   & =\left\{  N\;\Phi^{\prime}\left(  \det
x\right)  \det x-\Phi\left(  \det x\right)  \mid y_{j}^{n}=\Phi^{\prime
}\left(  \det x\right)  \det X_{n}^{k}\right\} \\
& =\left\{  N\;\Phi^{\prime}\left(  \psi\left(  \det y\right)  \right)
\psi\left(  \det y\right)  -\Phi\left(  \psi\left(  \det y\right)  \right)
\right\}
\end{align*}

yielding a unique value.
\end{proof}

\subsection{Example 1}

Take $\Phi\left(  t\right)  =\frac{N}{p}\left|  t\right|  ^{p/N}$, with $p\in
R$, so that
\begin{equation}
F\left(  x\right)  =\frac{N}{p}\left|  \det x\right|  ^{p/N}.\label{iv2}%
\end{equation}

Note that $F$ is homogeneous of degree $p$. If $p\neq0$, then $t^{N-1}%
\Phi^{\prime}\left(  t\right)  ^{N}=t^{p-1}$, which is invertible provided
$p\neq1$, yielding $\psi\left(  s\right)  =s^{1/\left(  p-1\right)  }$.
Substituting in the above, and taking advantage of the fact that $\Phi$ is
homogeneous of degree $p/N$, we get:
\begin{align*}
F^{L}(y)  & =(p-1)\Phi\left(  \psi\left(  \det y\right)  \right) \\
& =(p-1)\frac{N}{p}\left|  \det y\right|  ^{\frac{p}{p-1}\frac{1}{N}}%
\end{align*}

\begin{proposition}
If $p\neq0$ and $p\neq1$, the function $F:R^{N^{2}}\rightarrow R$ defined by
(\ref{iv2}) has a Legendre transform $F^{L}:R^{N^{2}}\rightarrow R$ defined
by:
\[
F^{L}(y)=\frac{N}{q}\left|  \det y\right|  ^{q/N},
\]
with $\frac{1}{p}+\frac{1}{q}=1$
\end{proposition}

Note that this duality holds for any value of $p$ different from $0$ and $1/N
$, including negative ones. Note also that if $N>1$ (the only interesting
case), this duality between $\Phi\left(  \det X\right)  $ and $\Psi\left(
\det Y\right)  $ has nothing to do with convexity. On the one hand, if
$\ p/N>1,$ so that $\Phi\left(  t\right)  $ is convex, then $0<q/N<1$, so that
$\Psi\left(  t\right)  $ is not convex. On the other hand, it is easy to check
that if $y\neq0$, we can find a matrix $\bar{x}$ such that $\sum x_{n}%
^{k}y_{k}^{n}=1$ and $\det\bar{x}=0$; it follows that the function
\[
x\rightarrow\sum x_{n}^{k}y_{k}^{n}-\frac{1}{p}\left(  \det X\right)  ^{p}%
\]
is unbounded from above and from below (consider the sequences $x_{n}=\pm
n\bar{x}$), and the critical point in the definition of the Legendre transform
cannot be a global minimum or maximum.

\subsection{Example 2}

Take $p=1$ in the above, so that $\Phi\left(  t\right)  =Nt^{1/N}$ and:
\begin{equation}
F\left(  x\right)  =N\left|  \det x\right|  ^{1/N}.\label{i8}%
\end{equation}

The function $F$ is defined on the whole of $R^{N^{2}}$. On the other hand, we
have
\[
y_{k}^{n}=\frac{\partial F}{\partial x_{n}^{k}}\left(  x\right)  =\Phi
^{\prime}\left(  \det x\right)  \det X_{n}^{k}%
\]
by (\ref{li7}), and $\det y=1$ by (\ref{li4}). So $F^{\prime}$ maps $R^{N^{2}%
}$ onto the set $\Sigma=\left\{  y\;|\;\det y=1\right\}  $. On the other hand,
formula (\ref{li11}) yields quite simply $\Gamma_{L}\left(  y\right)  =0$. Hence:

\begin{proposition}
The Legendre transform of the function $F:R^{N^{2}}\rightarrow R$ given by
(\ref{i8}) is the function $F^{L}:\Sigma\rightarrow R$ given by $F^{L}\left(
y\right)  \equiv0$
\end{proposition}

\subsection{Example 3}

Take $\Phi\left(  t\right)  =\ln\left|  t\right|  $, so that:
\begin{equation}
F\left(  x\right)  =\ln\left|  \det x\right| \label{ive}%
\end{equation}

Set $\Omega=\left\{  x\;|\;\det x\neq0\right\}  $.

\begin{proposition}
The function $F:\Omega\rightarrow R$ defined by (\ref{ive}) has a Legendre
transform $F^{L}:\Omega\rightarrow R$ defined by:
\[
F^{L}\left(  y\right)  =N+\ln\left|  \det y\right|
\]
\end{proposition}

The proof is left to the reader. It follows that the function $G\left(
x\right)  =\ln\left|  \det x\right|  +N/2$ is self-dual, i.e. $G^{L}=G$

\section{The case of $(N-1)$-cofactors}

We shall work with $N=K=3.$ We presume that similar results hold in the
general case, but we have not been able to handle the notations.

Denote by $x$the $3\times3$ matrix with coefficients $x_{n}^{k}$, with $1\leq
k\leq$ $K$ and $1\leq n\leq N$. Denote by $X_{n}^{k}$ the cofactor of
$x_{n}^{k}$ in $x$, and by $\Delta_{n}^{k}$ its determinant. Set
\[
\Delta=\left(  \Delta_{n}^{k}\right)  _{1\leq n\leq N}^{1\leq k\leq K}\in
R^{9}%
\]

Let $\Phi:R^{9}\rightarrow R$ be given. Consider the function:%

\begin{equation}
F\left(  x\right)  =\Phi(\Delta)\label{evu}%
\end{equation}

The formula for the Legendre transform then \ becomes:%

\begin{equation}
\Gamma_{L}\left(  y\right)  =\left\{  2\sum_{n,k}\frac{\partial\Phi}%
{\partial\Delta_{n}^{k}}\left(  \Delta\right)  \Delta_{n}^{k}-\Phi\left(
\Delta\right)  \mid y_{j}^{n}=\sum_{n,k}\frac{\partial\Phi}{\partial\Delta
_{n}^{k}}\frac{\partial\Delta_{n}^{k}}{\partial x_{n}^{j}}\right\} \label{Leg}%
\end{equation}

Let us simplify this formula a little bit by setting:
\[
\Phi_{k}^{n}=\frac{\partial\Phi}{\partial\Delta_{n}^{k}}\left(  \Delta\right)
\]

We then have:
\begin{equation}
y_{k}^{n}=\sum_{i,j}\Phi_{j}^{i}\frac{\partial\Delta_{i}^{j}}{\partial
x_{n}^{k}}=\sum_{\substack{i\neq n \\j\neq k}}\Phi_{j}^{i}\frac{\partial
\Delta_{i}^{j}}{\partial x_{n}^{k}}\left(  x\right) \label{cof}%
\end{equation}
because, if $k=j$ or $n=i$, the variable $x_{n}^{k}$ does not appear in the
cofactor $X_{n}^{k}$. If $k\neq j$, we shall denote by $p\left(  j,k\right)  $
the number in $\left\{  1,2,3\right\}  $ which is different from both $k$ and
$j$. Similarly, if $n\neq i$, we shall denote by $q\left(  i,n\right)  $ the
number in $\left\{  1,2,3\right\}  $ which is different from both $n$ and $i$.
If $k\neq j$ and $n\neq i$, we have:
\[
\Delta_{i}^{j}=\left(  -1\right)  ^{m\left(  i,j,k,n\right)  }\left(
x_{n}^{k}x_{q\left(  i,n\right)  }^{p\left(  j,k\right)  }-x_{q\left(
i,n\right)  }^{k}x_{n}^{p\left(  j,k\right)  }\right)
\]
where: the exponent $m\left(  i,j,k,n\right)  $ is $0$ if $k>p\left(
j,k\right)  $ and $n>q\left(  i,n\right)  $, or if $k<p\left(  j,k\right)  $
and $n<q\left(  i,n\right)  $, and $m\left(  i,j,k,n\right)  =$ $0$ otherwise.
It follows that:
\[
\frac{\partial\Delta_{i}^{j}}{\partial x_{n}^{k}}=\left(  -1\right)
^{m\left(  i,j,k,n\right)  }x_{q\left(  i,n\right)  }^{p\left(  j,k\right)  }%
\]
and hence:
\[
y_{k}^{n}=\sum_{\substack{i\neq n \\j\neq k}}\left(  -1\right)  ^{\left(
k-p\left(  j,k\right)  \right)  \left(  n-q\left(  i,n\right)  \right)  }%
\Phi_{j}^{i}x_{q\left(  i,n\right)  }^{p\left(  j,k\right)  }%
\]
\ 

Let us now consider the $3\times3$ matrix $y$ with coefficients $y_{k}^{n}$,
denote by $Y_{k}^{n}$ the cofactor of $y_{k}^{n}$ and by $D_{k}^{n}$ its
determinants$.$ Without loss of generality, we can assume that $n=k=1$, and we
get:
\begin{align*}
D_{1}^{1}  & =y_{2}^{2}y_{3}^{3}-y_{3}^{2}y_{2}^{3}\\
& =\left(  \Phi_{1}^{1}x_{3}^{3}-\Phi_{3}^{1}x_{3}^{1}-\Phi_{1}^{3}x_{1}%
^{3}+\Phi_{3}^{3}x_{1}^{1}\right)  \left(  \Phi_{1}^{1}x_{2}^{2}+\Phi_{2}%
^{1}x_{2}^{1}+\Phi_{1}^{2}x_{1}^{2}+\Phi_{2}^{2}x_{1}^{1}\right)  -\\
& \left(  -\Phi_{1}^{1}x_{3}^{2}-\Phi_{2}^{1}x_{3}^{1}+\Phi_{1}^{3}x_{1}%
^{2}+\Phi_{2}^{3}x_{1}^{1}\right)  \left(  -\Phi_{1}^{1}x_{2}^{3}+\Phi_{3}%
^{1}x_{2}^{1}-\Phi_{1}^{2}x_{1}^{3}+\Phi_{3}^{2}x_{1}^{1}\right) \\
& =\left(  \Phi_{1}^{1}\right)  ^{2}\left(  x_{3}^{3}x_{2}^{2}-x_{3}^{2}%
x_{2}^{3}\right)  +\left(  \Phi_{1}^{1}\Phi_{2}^{1}\right)  \left(  x_{3}%
^{3}x_{2}^{1}-x_{3}^{1}x_{2}^{3}\right)  +\left(  \Phi_{1}^{1}\Phi_{1}%
^{2}\right)  \left(  x_{3}^{3}x_{1}^{2}-x_{3}^{2}x_{1}^{3}\right)  +\\
& \left(  \Phi_{1}^{1}\Phi_{2}^{2}x_{3}^{3}x_{1}^{1}-\Phi_{2}^{1}\Phi_{1}%
^{2}x_{3}^{1}x_{1}^{3}\right)  +\left(  \Phi_{3}^{1}\Phi_{1}^{1}\right)
\left(  -x_{3}^{1}x_{2}^{2}+x_{3}^{2}x_{2}^{1}\right)  +\left(  -\Phi_{3}%
^{1}\Phi_{1}^{2}x_{3}^{1}x_{1}^{2}+\Phi_{1}^{1}\Phi_{3}^{2}x_{3}^{2}x_{1}%
^{1}\right)  +\\
& \left(  -\Phi_{3}^{1}\Phi_{2}^{2}x_{3}^{1}x_{1}^{1}+\Phi_{2}^{1}\Phi_{3}%
^{2}x_{3}^{1}x_{1}^{1}\right)  +\left(  \Phi_{1}^{3}\Phi_{1}^{1}\right)
\left(  -x_{1}^{3}x_{2}^{2}+x_{1}^{2}x_{2}^{3}\right)  +\left(  -\Phi_{1}%
^{3}\Phi_{2}^{1}x_{1}^{3}x_{2}^{1}+\Phi_{2}^{3}x_{1}^{1}\Phi_{1}^{1}x_{2}%
^{3}\right)  +\\
& \left(  -\Phi_{1}^{3}\Phi_{2}^{2}x_{1}^{3}x_{1}^{1}+\Phi_{2}^{3}x_{1}%
^{1}\Phi_{1}^{2}x_{1}^{3}\right)  +\left(  \Phi_{3}^{3}x_{1}^{1}\Phi_{1}%
^{1}x_{2}^{2}-\Phi_{1}^{3}x_{1}^{2}\Phi_{3}^{1}x_{2}^{1}\right)  +\left(
\Phi_{3}^{3}x_{1}^{1}\Phi_{2}^{1}x_{2}^{1}-\Phi_{2}^{3}x_{1}^{1}\Phi_{3}%
^{1}x_{2}^{1}\right)  +\\
& \left(  \Phi_{3}^{3}x_{1}^{1}\Phi_{1}^{2}x_{1}^{2}-\Phi_{1}^{3}x_{1}^{2}%
\Phi_{3}^{2}x_{1}^{1}\right)  +\left(  \Phi_{3}^{3}x_{1}^{1}\Phi_{2}^{2}%
x_{1}^{1}-\Phi_{2}^{3}x_{1}^{1}\Phi_{3}^{2}x_{1}^{1}\right)
\end{align*}

\begin{lemma}
If the matrix $\Phi_{k}^{n}$ has rank $1$, then the $D_{k}^{n}$ can be
expressed in terms of the $\Delta_{n}^{k}$ as follows:
\begin{equation}
D_{k}^{n}=\Phi_{k}^{n}\left(  \Phi_{3}^{2}\Delta_{2}^{3}+\Phi_{1}^{3}%
\Delta_{3}^{1}+\Phi_{2}^{3}\Delta_{3}^{2}+\Phi_{3}^{3}\Delta_{3}^{3}+\Phi
_{3}^{2}\Delta_{2}^{3}+\Phi_{1}^{3}\Delta_{3}^{1}+\Phi_{2}^{3}\Delta_{3}%
^{2}+\Phi_{3}^{3}\Delta_{3}^{3}\right) \label{hom}%
\end{equation}
\end{lemma}

\begin{proof}
If the matrix $\Phi_{k}^{n}$ has rank $1$, all its $2$-cofactors vanish, so
that $\Phi_{k}^{n}\Phi_{j}^{i}=\Phi_{j}^{n}\Phi_{k}^{i}$. The previous
expression then simplifies:
\begin{align*}
D_{1}^{1}  & =\left(  \Phi_{1}^{1}\right)  ^{2}\Delta_{1}^{1}+\left(  \Phi
_{1}^{1}\Phi_{2}^{1}\right)  \Delta_{1}^{2}+\left(  \Phi_{1}^{1}\Phi_{1}%
^{2}\right)  \Delta_{2}^{1}+\\
& \left(  \Phi_{1}^{1}\Phi_{2}^{2}\right)  \Delta_{2}^{2}+\left(  \Phi_{3}%
^{1}\Phi_{1}^{1}\right)  \Delta_{1}^{3}+\left(  \Phi_{1}^{1}\Phi_{3}%
^{2}\right)  \Delta_{2}^{3}+\\
& \left(  \Phi_{1}^{3}\Phi_{1}^{1}\right)  \Delta_{3}^{1}+\left(  \Phi_{2}%
^{3}\Phi_{1}^{1}\right)  \Delta_{3}^{2}+\left(  \Phi_{3}^{3}\Phi_{1}%
^{1}\right)  \Delta_{3}^{3}%
\end{align*}
and $\Phi_{1}^{1}$ factors out.
\end{proof}

If $\Phi$ is homogeneous of degree $\alpha$, the expression (\ref{hom})
simplifies by the Euler identity:
\begin{equation}
D_{k}^{n}=\alpha\Phi_{k}^{n}\Phi=\alpha\Phi\left(  \Delta\right)
\frac{\partial\Phi}{\partial\Delta_{n}^{k}}\left(  \Delta\right) \label{hon}%
\end{equation}
and the formula (\ref{Leg})for the Legendre transform $\Gamma_{L}$ of $F$ becomes:%

\[
\Gamma_{L}\left(  y\right)  =\left\{  \left(  2\alpha-1\right)  \Phi\left(
\Delta\right)  \mid D=\alpha\Phi\left(  \Delta\right)  \Phi^{\prime}\left(
\Delta\right)  \right\}
\]

\begin{proposition}
Assume that the function $F:R^{9}\rightarrow R$ is given by $F\left(
x\right)  =$ $\Phi\left(  \Delta\right)  $, where $\Phi$ is homogeneous of
degree $\alpha$, and the matrix $\partial\Phi/\partial\Delta_{n}^{k}$ has rank
$1$ everywhere. Assume that $\Sigma=F^{\prime}\left(  R^{9}\right)  $ is a
submanifold, and that:
\[
\left[  D_{1}=\alpha\Phi\left(  \Delta_{1}\right)  \Phi^{\prime}\left(
\Delta_{1}\right)  \text{ and }D_{2}=\alpha\Phi\left(  \Delta_{2}\right)
\Phi^{\prime}\left(  \Delta_{2}\right)  \right]  \Longrightarrow\Phi\left(
\Delta_{1}\right)  =\Phi\left(  \Delta_{2}\right)
\]
Then $F$ has a Legendre transform $F^{L}:\Sigma\rightarrow R$ given by
\begin{equation}
F^{L}\left(  y\right)  =\Psi\left(  D\right) \label{Lt}%
\end{equation}
where the $D_{k}^{n}$ are the determinants of the $2$-cofactors of $y$, and
$\Psi\left(  D\right)  =\left(  2\alpha-1\right)  \Phi\left(  \Delta\right)  $
for any $D$ such that $D=\alpha\Phi\left(  \Delta\right)  \Phi^{\prime}\left(
\Delta\right)  $
\end{proposition}

\subsection{Example 4}

We consider functions $F:R^{9}\rightarrow R$ of the following type:
\[
F\left(  x\right)  =\left(  (\sum_{n}\Delta_{n}^{1})^{\alpha}+(\sum_{n}%
\Delta_{n}^{2})^{\alpha}+(\sum_{n}\Delta_{n}^{3})^{\alpha}\right)  ^{\beta
}=\Phi\left(  \Delta\right)
\]

We have
\[
\Phi_{k}^{i}=\frac{\partial\Phi}{\partial\Delta_{i}^{k}}=\alpha\beta\left(
(\sum_{n}\Delta_{n}^{1})^{\alpha}+(\sum_{n}\Delta_{n}^{2})^{\alpha}+(\sum
_{n}\Delta_{n}^{3})^{\alpha}\right)  ^{\beta-1}(\sum_{n}\Delta_{n}%
^{k})^{\alpha-1}%
\]
so clearly the matrix $\Phi_{k}^{n}$ has rank $1$. The equations (\ref{hon})
become:
\[
D_{k}^{i}=\alpha\beta\Phi_{k}^{i}\Phi=\left(  \alpha\beta\right)  ^{2}\left(
(\sum_{n}\Delta_{n}^{1})^{\alpha}+(\sum_{n}\Delta_{n}^{2})^{\alpha}+(\sum
_{n}\Delta_{n}^{3})^{\alpha}\right)  ^{2\beta-1}(\sum_{n}\Delta_{n}%
^{k})^{\alpha-1}%
\]
from which we get
\begin{equation}
D_{k}^{1}=D_{k}^{2}=D_{k}^{3}\text{ for }k=1,2,3\label{77}%
\end{equation}

In other words, the Legendre transform will live on the $3$-dimensional
subspace $\Sigma$ of $R^{9}$ defined by the equations (\ref{77}). Setting
$D_{k}^{i}=D_{k}$ for every $i$, we continue the computations:%

\begin{align*}
\left(  D_{k}\right)  ^{\frac{\alpha}{\alpha-1}}  & =\left(  \alpha
\beta\right)  ^{\frac{2\alpha}{\alpha-1}}\left(  (\sum_{n}\Delta_{n}%
^{1})^{\alpha}+(\sum_{n}\Delta_{n}^{2})^{\alpha}+(\sum_{n}\Delta_{n}%
^{3})^{\alpha}\right)  ^{\left(  2\beta-1\right)  \frac{\alpha}{\alpha-1}%
}(\sum_{n}\Delta_{n}^{k})^{\alpha}\\
\sum_{k}\left(  D_{k}\right)  ^{\frac{\alpha}{\alpha-1}}  & =\left(
\alpha\beta\right)  ^{\frac{2\alpha}{\alpha-1}}\left(  (\sum_{n}\Delta_{n}%
^{1})^{\alpha}+(\sum_{n}\Delta_{n}^{2})^{\alpha}+(\sum_{n}\Delta_{n}%
^{3})^{\alpha}\right)  ^{\frac{2\alpha\beta-1}{\alpha-1}}\\
\left(  \sum_{k}\left(  D_{k}\right)  ^{\frac{\alpha}{\alpha-1}}\right)
^{\frac{\beta\left(  \alpha-1\right)  }{2\alpha\beta-1}}  & =\left(
\alpha\beta\right)  ^{\frac{2\alpha\beta}{2\alpha\beta-1}}\left(  (\sum
_{n}\Delta_{n}^{1})^{\alpha}+(\sum_{n}\Delta_{n}^{2})^{\alpha}+(\sum_{n}%
\Delta_{n}^{3})^{\alpha}\right)  ^{\beta}%
\end{align*}

Finally, the Legendre transform of $F$ turns out to be the function:
\[
F^{L}\left(  y\right)  =(2\alpha\beta-1)\left(  \alpha\beta\right)
^{-\frac{2\alpha\beta}{2\alpha\beta-1}}\left(  \sum_{k}\left(  D_{k}\right)
^{\frac{\alpha}{\alpha-1}}\right)  ^{\frac{\beta\left(  \alpha-1\right)
}{2\alpha\beta-1}}%
\]
restricted to the $3$-dimensional subspace $\Sigma\subset R^{9}$ defined by
the relations $D_{k}^{n}=D_{k}.$ Here, $D_{k}^{n}$ denotes the cofactor of
$y_{k}^{n}$ in the matrix $Y$.

Note that $\ F$ is homogeneous of degree $2\alpha\beta$ and $F^{L}$ is
homogeneous of degree $2\alpha\beta/\left(  2\alpha\beta-1\right)  $. Setting
$p=2\alpha\beta$ and $q=2\alpha\beta/\left(  2\alpha\beta-1\right)  ,$ we find
that;
\[
\frac{1}{p}+\frac{1}{q}=1
\]

\subsection{Example 5}

Let $\left(  x^{1},x^{2}\right)  $ be a pair of vectors in $R^{3}$. We
consider functions $F:R^{6}\rightarrow R$ of the following type:
\[
F\left(  x^{1},x^{2}\right)  =\Phi\left(  \det\left|
\begin{array}
[c]{cc}%
x_{1}^{1} & x_{1}^{2}\\
x_{2}^{1} & x_{2}^{2}%
\end{array}
\right|  ,\det\left|
\begin{array}
[c]{cc}%
x_{1}^{1} & x_{1}^{2}\\
x_{3}^{1} & x_{3}^{2}%
\end{array}
\right|  ,\det\left|
\begin{array}
[c]{cc}%
x_{2}^{1} & x_{2}^{2}\\
x_{3}^{1} & x_{3}^{2}%
\end{array}
\right|  \right)
\]
where $\Psi:R^{3}\rightarrow R$ is a $C^{1}$ function. In the previous
framework, this can be understood as a function $F\left(  X\right)  $, where
$X $ is a $3\times3$ matrix, which depends only on the first three cofactors.
Clearly the rank condition will hold, and the previous results apply. It will
be more convenient, however, to run through the computations again in that
particular case, with simplified notations.

Set $\Delta=\left(  \Delta_{3},\Delta_{2},\Delta_{1}\right)  ,$ with:
\begin{align*}
\Delta_{3}  & =x_{1}^{1}x_{2}^{2}-x_{1}^{2}x_{2}^{1}\\
\Delta_{2}  & =x_{1}^{1}x_{3}^{2}-x_{1}^{2}x_{3}^{1}\\
\Delta_{1}  & =x_{2}^{1}x_{3}^{2}-x_{2}^{2}x_{3}^{1}%
\end{align*}

We have $F\left(  x^{1},x^{2}\right)  =\Phi\left(  \Delta\right)  $. Set
$y_{k}^{n}=\partial F/\partial x_{n}^{k}$ and compute the cofactors. We get:%

\begin{align*}
D^{3}  & =y_{1}^{1}y_{2}^{2}-y_{2}^{1}y_{1}^{2}=\left(  \frac{\partial\Phi
}{\partial\Delta_{3}}\right)  ^{2}\Delta_{3}+\left(  \frac{\partial\Phi
}{\partial\Delta_{3}}\frac{\partial\Phi}{\partial\Delta_{1}}\right)
\Delta_{1}+\left(  \frac{\partial\Phi}{\partial\Delta_{2}}\frac{\partial\Phi
}{\partial\Delta_{3}}\right)  \Delta_{2}\\
D^{2}  & =y_{1}^{1}y_{2}^{3}-y_{2}^{1}y_{1}^{3}=\left(  \frac{\partial\Phi
}{\partial\Delta_{2}}\right)  ^{2}\Delta_{2}+\left(  \frac{\partial\Phi
}{\partial\Delta_{3}}\frac{\partial\Phi}{\partial\Delta_{2}}\right)
\Delta_{3}+\left(  \frac{\partial\Phi}{\partial\Delta_{2}}\frac{\partial\Phi
}{\partial\Delta_{1}}\right)  \Delta_{1}\\
D^{1}  & =y_{1}^{2}y_{2}^{3}-y_{2}^{2}y_{1}^{3}=\left(  \frac{\partial\Phi
}{\partial\Delta_{1}}\right)  ^{2}\Delta_{1}+\left(  \frac{\partial\Phi
}{\partial\Delta_{1}}\frac{\partial\Phi}{\partial\Delta_{2}}\right)
\Delta_{2}+\left(  \frac{\partial\Phi}{\partial\Delta_{1}}\frac{\partial\Phi
}{\partial\Delta_{3}}\right)  \Delta_{3}%
\end{align*}

We summarize:
\begin{equation}
D^{n}=\frac{\partial\Phi}{\partial\Delta_{n}}\left[  \frac{\partial\Phi
}{\partial\Delta_{1}}\Delta_{1}+\frac{\partial\Phi}{\partial\Delta_{2}}%
\Delta_{2}+\frac{\partial\Phi}{\partial\Delta_{3}}\Delta_{3}\right]
,n=1,2,3\label{s0}%
\end{equation}

As a particular case, consider the function:
\begin{equation}
F\left(  x^{1},x^{2}\right)  =\left[  \left(  x_{1}^{1}x_{2}^{2}-x_{1}%
^{2}x_{2}^{1}\right)  ^{\alpha}+\left(  x_{1}^{1}x_{3}^{2}-x_{1}^{2}x_{3}%
^{1}\right)  ^{\alpha}+\left(  x_{2}^{1}x_{3}^{2}-x_{2}^{2}x_{3}^{1}\right)
^{\alpha}\right]  ^{1/\beta}\label{area}%
\end{equation}

When $\alpha=1/\beta=2$, this gives the area of the triangle spanned by the
vectors $x^{1}$ and $x^{2}$. We apply the preceding result, with:
\[
\Phi\left(  \Delta_{3},\Delta_{2},\Delta_{1}\right)  =\left[  \left(
\Delta_{1}\right)  ^{\alpha}+\left(  \Delta_{2}\right)  ^{\alpha}+\left(
\Delta_{3}\right)  ^{\alpha}\right]  ^{\beta}%
\]

The system\ (\ref{s0}) becomes:%

\[
D^{n}=(\alpha\beta)^{2}\left[  \left(  \Delta_{1}\right)  ^{\alpha}+\left(
\Delta_{2}\right)  ^{\alpha}+\left(  \Delta_{3}\right)  ^{\alpha}\right]
^{2\beta-1}\left(  \Delta_{n}\right)  ^{\alpha-1}%
\]
and can easily be inverted (note that if $\alpha=1/\beta=2$, we get the
identity). We get:%

\[
\Delta_{n}=(\alpha\beta)^{-\frac{2\alpha-1}{2\alpha\beta-1}\frac{1}{\alpha-1}%
}\left[  \left(  D^{1}\right)  ^{\frac{\alpha}{\alpha-1}}+\left(
D^{2}\right)  ^{\frac{\alpha}{\alpha-1}}+\left(  D^{3}\right)  ^{\frac{\alpha
}{\alpha-1}}\right]  ^{-\frac{2\beta-1}{2\alpha\beta-1}}\left(  D^{n}\right)
^{\frac{1}{\alpha-1}}%
\]
Substituting into formula (\ref{Lt}), and taking advantage of the fact that
$\Phi$ is homogeneous of degree $\alpha\beta$, we get the Legendre transform:
\begin{align*}
F^{L}\left(  y_{1},y_{2}\right)   & =(2\alpha\beta-1)\left[  \left(
\Delta_{1}\right)  ^{\alpha}+\left(  \Delta_{2}\right)  ^{\alpha}+\left(
\Delta_{3}\right)  ^{\alpha}\right]  ^{\beta}\\
& =(2\alpha\beta-1)(\alpha\beta)^{-\frac{2\alpha-1}{2\alpha\beta-1}%
\frac{\alpha\beta}{\alpha-1}}\left[  \left(  D^{1}\right)  ^{\frac{\alpha
}{\alpha-1}}+\left(  D^{2}\right)  ^{\frac{\alpha}{\alpha-1}}+\left(
D^{3}\right)  ^{\frac{\alpha}{\alpha-1}}\right]  ^{\frac{\alpha-1}%
{2\alpha\beta-1}\beta}\\
& =(2\alpha\beta-1)(\alpha\beta)^{-\frac{2\alpha-1}{2\alpha\beta-1}%
\frac{\alpha\beta}{\alpha-1}}\left[  \left(  y_{1}^{2}y_{2}^{3}-y_{2}^{2}%
y_{1}^{3}\right)  ^{\frac{\alpha}{\alpha-1}}+\left(  y_{1}^{1}y_{2}^{3}%
-y_{2}^{1}y_{1}^{3}\right)  ^{\frac{\alpha}{\alpha-1}}+\left(  y_{1}^{1}%
y_{2}^{2}-y_{2}^{1}y_{1}^{2}\right)  ^{\frac{\alpha}{\alpha-1}}\right]
^{\frac{\alpha-1}{2\alpha\beta-1}\beta}%
\end{align*}

Note that if $\alpha=1/\beta=2$, we find $F=F^{L}$: the function $F$ is its
own Legendre transform. Note also that $F$ is homogeneous of degree
$2\alpha\beta$ and $F^{L}$ homogeneous of degree $2\alpha\beta/\left(
2\alpha\beta-1\right)  .$ Setting $p=2\alpha\beta$ and $q=2\alpha\beta/\left(
2\alpha\beta-1\right)  ,$ we find that;
\[
\frac{1}{p}+\frac{1}{q}=1
\]
as before.

\section{A variational problem.}

As a example of possible application of this kind of duality, let us consider
the following problem. Given a positive definite quadratic form $\left(
Ax,x\right)  $ on $R^{N^{2}}$, and a $N\times N$-matrix $f,$ we want to solve
$\Phi^{\prime}\left(  x\right)  =0$, where::
\[
\Phi\left(  x\right)  =\frac{1}{2}\left(  Ax,x\right)  -\frac{N}{p}\left|
\det x\right|  ^{p/N}-\left(  f,x\right)
\]

Such points are called critical points of $\Phi$. Any critical point of $F$
solves the system:
\[
Ax=\left|  \det x\right|  ^{p/N-1}X+f
\]
where $X$ is the matrice of cofactors of $x$.

\begin{proposition}
\label{p11}If $p<2$, there is at least one critical point for $\Phi$.
\end{proposition}

\begin{proof}
Since $p<2$, the function $\Phi$ is coercive:\ $\Phi\left(  x\right)
\rightarrow\infty$ when $\left\|  x\right\|  \rightarrow\infty$. So it attains
its minimum at some $x$, which has to be a critical point.
\end{proof}

We now use duality theory to treat the case $p>2.$ We shall use an extension
of the Clarke duality formula (see \cite{ie} ):

\begin{proposition}
\label{p12}Suppose $F\left(  x\right)  $ has a Legendre transform
$F^{L}\left(  y\right)  $. Consider the functions $\Phi$ and $\Psi$ defined
by:
\begin{align*}
\Phi\left(  x\right)   & =\frac{1}{2}\left(  Ax-f,x\right)  -F\left(  x\right)
\\
\Psi\left(  y\right)   & =\frac{1}{2}\left(  Ay+f,y\right)  -F^{L}\left(
Ay\right)
\end{align*}
If $y$ is a critical point of $\Psi$, then $x=y+A^{-1}f$ is a critical point
of $\Phi.$
\end{proposition}

\begin{proof}
If $y$ is a critical point of $\Psi$, we have $\Psi^{\prime}\left(  y\right)
=0$, and hence $Ay+f=A\left(  F^{L}\right)  ^{\prime}\left(  Ay\right)  $.
Since $A$ is invertible, it follows that $y+A^{-1}f=\left(  F^{L}\right)
^{\prime}\left(  Ay\right)  $. Since $\left[  \left(  F^{L}\right)  ^{\prime
}\right]  ^{-1}=F^{\prime}$, it follows that
\[
F^{\prime}\left(  y+A^{-1}f\right)  =Ay=A\left(  y+A^{-1}f\right)  -f
\]
so that $F^{\prime}\left(  x\right)  =Ax-f$, as desired
\end{proof}

\begin{proposition}
If $p>2$, there is at least one non-trivial critical point for $\Phi$.
\end{proposition}

\begin{proof}
Consider the function:
\[
\Psi\left(  y\right)  =\frac{1}{2}\left(  Ay+f,y\right)  -\frac{N}{q}\left|
\det Ay\right|  ^{q/N}%
\]
By proposition \ref{p11}, the function $\Psi$ has a critical point $\bar
{y}\neq0.$ By proposition \ref{p12} it is also a critical point of $\Phi$
\end{proof}

Note that in this case $\inf\Phi=-\infty$, so that the critical point $x$
cannot be a minimizer.

\end{document}